\documentclass[a4paper]{article}
\usepackage[dvips]{graphics}
\usepackage{mathrsfs}
\usepackage{amssymb}
\usepackage{stmaryrd}
\usepackage{pifont}
\usepackage{amsmath}
\usepackage{hyperref}
\usepackage{amsthm}
\usepackage{cases}
\allowdisplaybreaks 
\def\beqn{\begin{displaymath}}
\def\eeqn{\end{displaymath}}

\def\bea{\begin{eqnarray}}
\def\eea{\end{eqnarray}}
\def\bean{\begin{eqnarray*}}
\def\eean{\end{eqnarray*}}

\def\beq{\begin{equation}}
\def\eeq{\end{equation}}

\begin{document}

\begin{center}
{\Large   The Stability of the Solutions for a Quasilinear Degenerate Parabolic
 Equation }
\bigskip

Miao Ouyang

School of Applied Mathematics, Xiamen University of Technology,
Xiamen 361024, China\ \

Email:mouyang@xmut.edu.cn

\ Huashui Zhan

School of Applied Mathematics, Xiamen University of Technology,
Xiamen 361024, China\ \

Email: huashuizhan@163.com

\end{center}

\textbf{Abstract} \ The equation arising from Prandtl boundary layer theory
 $$
\frac{\partial u}{\partial t} -\frac{\partial }{\partial x_i}\left( a(u,x,t)\frac{\partial u}{\partial x_i}\right)
-f_i(x)D_iu+c(x,t)u=g(x,t)
$$
is considered. The existence of the entropy solution can be proved by BV estimate method. The interesting problem is that,  since $a(\cdot,x,t)$ may be degenerate on the boundary,  the usual boundary value condition may be overdetermined. Accordingly, only dependent on a partial boundary value condition, the stability of solutions can be expected. This expectation is turned to reality by Kru\u{z}kov's bi-variables method, a reasonable partial boundary value condition matching up with the equation is found first time.  Moreover, if $a_{x_i}(\cdot,x,t)\mid_{x\in \partial \Omega}=a(\cdot,x,t)\mid_{x\in \partial \Omega}=0$ and $f_i(x)\mid_{x\in \partial \Omega}=0$, the stability can be proved even without any boundary value condition.


 \textbf{2000 MR Subject Classification} 35K65, 35L65,  35R35

\bigskip
\section{Introduction }
\bigskip

The initial-boundary value problem of the quasilinear degenerate parabolic equation
$$
\frac{\partial u}{\partial t} -\frac{\partial }{\partial x_i}\left( a(u,x,t)\frac{\partial u}{\partial x_i}\right)
-f_i(x)D_iu+c(x,t)u=g(x,t),\ \ (x,t)\in \Omega \times (0,T),\eqno(1.1)
$$
 \begin{equation*}
u(x,0)=u_0(x), \ x\in \Omega,\eqno(1.2)
\end{equation*}
$$
u(x,t)=0, (x,t)\in \partial\Omega.\eqno(1.3)
$$
is considered in this paper, where $a(u,x,t)\geq 0$, $\Omega \subset \mathbb{R}^{N}$ is a appropriately smooth open domain, $D_i =\frac{\partial}{\partial x_i}$, the double indices of $i$  represent the summation from $1$ to $N$ as usual.

Equation (1.1) arises from the boundary layer theory [1] etc. As the simplification of the
Navier-Stokes equation, the Prandtl boundary layer equation describes the motion of a fluid with small
viscosity about a solid body in a thin layer which is formed near
its surface owing to the adhesion of the viscous fluid to the solid
surface. In particular, we consider the motion of a fluid occupying a two
dimensional region is characterized by the velocity vector
$V=(u,v)$, where $u,v$ are the projections of $V$ onto the
coordinate axes $x,y,$ respectively, assume that the
density of the fluid $\rho$ is equal to 1, then the Prandtl boundary layer equation for a
non-stationary boundary layer arising in an axially symmetric
incompressible flow past a solid body has the form [1]
\begin{equation*}
\left\{\begin{array}{ll}
u_{t}+uu_{x}+vu_{y}=\nu u_{yy}-p_{x}               \mbox{,}\\
u_{x}+v_{y}=0\mbox{,}\\
u(0,x,y)=u_{0}(x,y),\quad u(t,0,y)=u_{1}(t,y){,}\\
 u(t,x,0)=0,\quad v(t,x,0)=v_{0}(t,x){,}\\
 \lim_{y \to \infty} u(t,x,y)=U(t,x){.}
\end{array}\right.
\end{equation*}
in a domain $D=\{0<t<T,0<x<X,0<y<\infty\}$, where $\nu=const>0$ is
the viscosity coefficient of the incompressible fluid, $u_{0}>0,u_{1}>0$
 for $y>0$, $u_{0y}>0, u_{1y}>0$ for $y\geq0$, where, $p=p(t,x)$ is the pressure, $U=U(t,x)$ is the velocity at the outer edge of the boundary layer which satisfies
 $$
 U_{t}+UU_{x}=-p_{x}(t,x),\  U(t,x)>0.
 $$
By the well-known Crocco transform,
\begin{equation*}
\tau=t,\quad\xi=x,\quad\eta=u(t,x,y),\quad w(\tau,\xi,\eta)=u_{y}.
\end{equation*}
 we can show that $u_{y}=w$ satisfies the following nonlinear equation
\begin{equation*}
w_{\tau}=\nu w^{2}w_{\eta\eta}-\eta
w_{\xi}+p_{x}w_{\eta}.\eqno(1.4)
\end{equation*}
By a linearized method, Oleinik had shown that there is a local classical solution to this equation [2]. Although there are some important papers to studied the global solutions of the Prandtl boundary layer equation [32-37], the related problems are far from being solved.  For example,  the compatibility problem between Navier-Stokes equation and Prandtl boundary layer equation. For another example,  whether there is a global solution of equation (1.4) and whether this global solution can be deduced a global weak solution of the Prandtl boundary layer equation by the inverse transform of Crocco transform ? In fact, if the domain is not the $N-$dimmensional cube, whether the inverse transform of Crocco transform exists or not is still unsolved. In addition, many reaction-diffusion problems can be summed up to equation (1.1) [2].

In this paper, we will consider the global solutions of equation (1.1). After the pioneering work [3] by Vol$^{\prime }$pert-Hudjaev,  the Cauchy problem of equation (1.1)  had been studied in [4-13] etc.,  the solutions to the Cauchy problem of equation (1.1) are well-posedness. Also, the initial-boundary value problem of equation (1.1) had been studied in many papers, many excellent and important results had been obtained in [14-16], [30-31] etc. Shall we say, there is not important problem left? I think it is too early to make such a conclusion. Besides the problems related to Prandtl boundary layer theory, since $a(u,x,t)\geq 0$ and may be degenerate in the interior of $\Omega$ or on the boundary $\partial \Omega$, everyone knows that the boundary value condition (1.3) is overdetermined, there is not an effective method to find a reasonable partial boundary value condition
\begin{equation*}
u(x,t)=0, (x,t)\in \Sigma_p \times (0,T),\eqno(1.5)
\end{equation*}
to replace (1.3), where $\Sigma_p$ is a relative open subset of $\partial \Omega$. Here, we like to suggest that the boundary value condition (1.3) or (1.5) is understood in the sense of the trace, and  we expect to find a analytic expression of $\Sigma_p$ in this paper. The difficulty comes from that, since the equation has the nonlinearity, the partial boundary $\Sigma_p$ in (1.5) can not be depicted out by Fichera function as that of the linear degenerate parabolic equation [26-27].

In fact, for a nonlinear parabolic equation, how to impose a reasonable partial boundary value condition has been up in the air for a long time [14-19]. Let us give some details. In [14-16], the entropy solutions defined in these references are in $L^{\infty}(Q_T)$ sense, one can not define the trace on the boundary, accordingly, it is impossible to express $\Sigma_p$ in an analytic formula. Instead, the authors of [14-16] had found a kind of the entropy inequality to imply the boundary value condition (1.5) in  ingenious ways. In the work by Yin-Wang [17], the degenerate non-Newtonian fluid equation
\begin{equation*}
\frac{\partial u}{\partial t}-\text{div}(d(x)|\nabla u|^{p-2}\nabla u)-f_i(x)D_iu+c(x,t)u=g(x,t),(x,t)\in \Omega\times (0,T),\eqno(1.6)
\end{equation*}
was considered.  By means of a reasonable integral description, in [17], the boundary $\partial \Omega$ is classified  into three parts: the nondegenerate boundary $\Sigma_1$, the weakly degenerate boundary $\Sigma_2$ and the strongly degenerate boundary $\Sigma_3$. Instead of the usual boundary condition (1.3),
a partial boundary value condition (1.5) is imposed, where
$$
\Sigma_p =\Sigma_2\bigcup\Sigma_3. \eqno(1.7)
$$
 It is pity that, since equation (1.1) is apparently different the Non-Newtonian equation (1.6), $\Sigma_p$ also can not be described as (1.7). If the domain $\Omega$ is the $N-$dimensional cube or the half space of $\mathbb{R}^N$,  the equation
 $$
\frac{\partial u}{\partial t} =\Delta A(u)
+\text{div}(b(u)),\ \ (x,t)\in \Omega \times (0,T),
$$
was studied in [18-19] by the author recently,  a reasonable analytic expression of $\Sigma_p$ had been found in [18-19]. However, for a general domain $\Omega$, the problem remains open. We hope to make a essential progress sooner or later.

Certainly, since the subset set $D_0=\{x\in \Omega: a(\cdot,x,t)=0\}$ may have a positive measure in $\Omega$, equation (1.1) has  hyperbolic characteristic in $D_0$. Thus, only in the sense of the entropy solution, the uniqueness (or the stability) of the weak solution can be obtained [1]. In this paper, with the help of the entropy solutions defined in the sense of $BV$ functions [1, 5, 18, 22], we study the well-posed problem of equation (1.1) with the initial value  (1.2) and the partial boundary value condition (1.5), the key is to find a reasonable analytic expression of $\Sigma_p$ first time.

   The paper is arranged as follows. After the introduction section, section 2 introduces the definition of the entropy solution and the main results. Section 3 gives the proof of the existence of the entropy solutions. Section 4 introduces the well-known Kru\u{z}kov bi-variables method. Section 5 is on the stability of the entropy solutions based on the partial boundary value condition. At the end, an explanation of the definition of the entropy solution is given.

  \bigskip
\section{The definition of the entropy solution and the main results}
\bigskip

For the completeness of the paper, we first quote the definition of $BV$ function and its properties [28].

\textbf{Definition 2.1} Let $\Omega\subset \mathbb{R}^{m}$ be an
open set and let $f\in L^{1}(\Omega)$. Define
$$
\int_{\Omega}|Df|=\sup\left\{\int_{\Omega}f \text{div} gdx:
g=(g_{1}, g_{2}, \cdots, g_{N})\in C^{1}_{0}(\Omega;
\mathbb{R}^{m}), |g(x)|\leq 1, x\in \Omega\right\},
$$
where $\text{div} g=\sum_{i=1}^{m}\frac{\partial g_{i}}{\partial
x_{i}}$.

\textbf{Definition 2.2} A function of $f\in L^{1}(\Omega)$ is said
to have bounded variation in $\Omega$ if
$$
\int_{\Omega}|Df|<\infty.
$$
We define $BV(\Omega)$ as the space of all functions in
$L^{1}(\Omega)$ with bounded variation.

This is equivalent to that the generalized derivatives of every
function in $ BV(\Omega)$ are regular measures on $\Omega$. Under
the norm
$$
\|f\|_{BV}=\|f\|_{L^{1}}+\int_{\Omega}|Df|,
$$
$BV(\Omega)$ is a Banach space.

\textbf{Proposition 2.3} (Semicontinuity)\  Let $\Omega\subseteq
\mathbb{R}^{m}$ be an open set and $\{f_{j}\}$ a sequence of
functions in $BV(\Omega)$ which converge in
$L^{1}_{\text{loc}}(\Omega)$ to a function $f$. Then
$$
\int_{\Omega}|Df|\leq\lim_{j\rightarrow \infty}\inf
\int_{\Omega}|Df_{j}|.
$$

\textbf{Proposition 2.4} (Integration by part)\  Let
$$
C^{+}_{R}=\mathscr{B}(0,R)\times(0,R)=\mathscr{B}_{R}\times(0,R)
$$
and $ f\in BV(C^{+}_{R})$. Then there exists a function $f^{+}\in
L^{1}(\mathscr{B}_{R})$ such that  for $H_{n-1}$-almost all $y\in
\mathscr{B}_{R}$,
$$
\lim_{\rho\rightarrow
0}\rho^{-m}\int_{C^{+}_{\rho}(y)}|f(z)-f^{+}(y)|dz=0.
$$
Moreover, if $C_{R}=\mathscr{B}_{R}\times(-R, R)$, then for every
$g\in C_{0}^{1}(C_{R}; \mathbb{R}^{m})$,
$$
\int_{C^{+}_{R}}f\text{div} gdx=-\int_{C^{+}_{R}}\langle g,
Df\rangle+\int_{\mathscr{B}_{R}}f^{+}gdH_{n-1},
$$
where $\mathscr{B}_{\rho }=\{x\in \mathbb{R}^{m};\mid x\mid <\rho
\}$.

\textbf{Remark 2.5} The function $f^{+}$ is called the trace of $f$
on $\mathscr{B}_{R}$ and obviously
$$
f^{+}(y)=\lim_{\rho\rightarrow
0}\frac{1}{|C^{+}_{\rho}(y)|}\int_{C^{+}_{\rho}(y)}f(z)dz.
$$
The definition of the trace is easy generalized to a general smooth domain in  $\mathbb{R}^{m}$.

Secondly, we give the definition of the entropy solutions matching up with equation (1.1). For small $\eta>0$, let
$$
S_{\eta}(s)=\int_{0}^{s}h_{\eta}(\tau)d\tau,
\ h_{\eta}(s)=\frac{2}{\eta}\left(1-\frac{\mid s\mid}{\eta
}\right)_{+}.
$$
Obviously $h_{\eta}(s)\in C(\mathbb{R})$, and
\begin{equation*}
h_{\eta}(s)\geq 0,\;\ \mid sh_{\eta}(s)\mid \leq 1,\;\mid S_{\eta
}(s)\mid \leq 1;\underset{\eta \rightarrow 0}{\lim}S_{\eta
}(s)=\text{sgn}s,\;\underset{\eta \rightarrow
0}{\lim}sS_{\eta}^{\prime}(s)=0.\eqno(2.1)
\end{equation*}

\textbf{Definition 2.6} A function $u\in BV(Q_{T})\cap L^{\infty}(Q_{T})$ is said to be the entropy
solution of equation (1.1) with the initial value (1.2), provided that

\ \ 1. There exist $g^i\in L^2(Q_T)$ $(i=1,2,\cdots, N)$ such that for any $\varphi(x,t)\in C^1_{0}(Q_T)$
$$
\iint\nolimits_{Q_{T}}\varphi(x,t)g^i(x,t)dxdt=\int\int\nolimits_{Q_{T}}\varphi(x,t)\widehat{\sqrt{a(u,x,t)}}\frac{\partial u}{\partial x_i}dxdt,\eqno(2.2)
$$
where
$$
\widehat{\sqrt{a(u,x,t)}}(u,x,t)=\int_{0}^{1}\sqrt{a(\tau u^{+}+(1-\tau )u^{-},x,t)}d\tau,
$$
 is the composite mean value of $\sqrt{a(u,x,t)}$.

\ \ 2. If $\varphi\in C_{0}^{2}(Q_{T})$ and
$\varphi\geq 0$, for $k\in \mathbb{R}$ and for any small $\eta>0$ there holds
$$
\begin{gathered}
\iint\nolimits_{Q_{T}}\left[I_{\eta}(u-k)\varphi_{t}-f_i(x)I_{\eta}(u-k)\varphi_{x_{i}}+A_{\eta}(u,x,t,k)\Delta \varphi
-\sum_{i=1}^{N}S_{\eta}^{\prime}(u-k)\mid g^{i}\mid
^{2}\varphi\right]dxdt
\hfill \\
 \begin{array}{*{20}{c}}
{}
\end{array}
+\iint \nolimits_{Q_{T}}\int_{k}^{u}a_{x_i}(s,x,t)S_{\eta}(s-k)ds\varphi_{x_{i}}dxdt
\hfill \\
 \begin{array}{*{20}{c}}
{}
\end{array}
-\iint\nolimits_{Q_{T}}f_{ix_i}(x)(u-k)\varphi S_{\eta
}(u-k)dxdt+\iint\nolimits_{Q_{T}}f_{ix_i}(x)\int_{k}^{u}(s-k)h_{\eta}(s-k)ds\varphi dxdt
\hfill \\
 \begin{array}{*{20}{c}}
{}
\end{array}
-\iint \nolimits_{Q_{T}}[c(x,t)u+g(x,t)]\varphi S_{\eta}(u-k) dxdt
\hfill \\
 \begin{array}{*{20}{c}}
{}
\end{array}
\geq 0.
\hfill\\
\end{gathered}
\eqno(2.3)
$$

\ \ 3. The initial value is satisfied in the sense of that
\begin{equation*}
\lim_{t\rightarrow 0}\int_{\Omega }\mid
u(x,t)-u_{0}(x)\mid dx=0.\eqno(2.4)
\end{equation*}
Here $\vec{f}=\{f_i\}$,
\begin{eqnarray*}
A_{\eta}(u,x,t, k) &=& \int_{k}^{u}a(s,x,t)S_{\eta}(s-k)ds,
\end{eqnarray*}
and
$$
 I_{\eta}(u-k)=\int_{0}^{u-k}S_{\eta}(s)ds.
$$

\textbf{Definition 2.7} If  $u\in BV(Q_{T})\cap L^{\infty}(Q_{T})$ is the entropy
solution of equation (1.1) with the initial value (1.2),  and the partial boundary value condition (1.5) is satisfied in the sense of the trace, then we say $u$ is a entropy solution of the initial-boundary value problem of equation (1.1). Here,
$$
\Sigma_p=\left\{x\in \partial \Omega: \sum_{i=1}^Nf_i(x)n_i<0\right\}\bigcup \left\{x\in \partial \Omega: \sum_{i=1}^Na_{x_i}(\cdot,x,t)n_i\neq 0\right\} \bigcup \left\{x\in \partial \Omega: a(\cdot,x,t)\neq 0\right\}, \eqno(2.5)
$$
and $\vec{n}=\{n_i\}$ is the inner normal vector of $\Omega$.

 In what follows, we can show that if $a(\cdot,x,t)\mid_{x\in\partial \Omega}=0$,  then $\Sigma_p$ in the partial boundary value (1.5) can be depicted out as (2.5).  Based on this fact, thirdly, we will prove the following theorems.

\textbf{Theorem 2.8} If $a(s,x,t) \in C^{1}(\mathbb{R}^N\times Q_T)$, $f_{i}(x) \in C(\overline{\Omega})$, $c(x,t)$ and $g(x,t)$ are $C^{1}(\overline{Q_T})$, $u_{0}(x)\in L^{\infty}(\Omega)$,
then equation (1.1) with the initial value condition
(1.2) has an entropy solution in the sense of Definition 2.6.

\textbf{Theorem 2.9} If $a(s,x,t) \in C^{1}(\mathbb{R}^N\times Q_T)$ with $a(0,x,t) =0, (x,t)\in Q_{T}$, $f_{i}(x) \in C^{1}(\overline{\Omega})$, $c(x,t)$ and $g(x,t)$ are $C^{1}(\overline{Q_T})$, $u_{0}(x)\in L^{\infty}(\Omega)$,
and there is a constant $\delta_1>0$ such that
\begin{equation*}
a(r,x,t)-\delta_1\sum_{s=1}^{N+1}(a_{x_{s}}(r,x,t))^{2}\geq 0,\eqno(2.6)
\end{equation*}
then the initial-boundary value problem of equation (1.1) has an entropy solution in the sense of Definition 2.7.

\textbf{Theorem 2.10} Suppose $a(\cdot,x,t)$ is a $C^{1}(\overline{Q_T})$ function, $a(s,x,t)$ is bounded when $s$ is bounded,  $f_{i}(x) \in C^{1}(\overline{\Omega})$,  $c(x,t)$ and $g(x,t)$ are bounded. Suppose that when $x$ is near to the boundary,
$$
\Delta d\leq 0,\eqno(2.7)
$$
there exist constants $\delta_2>00$ such that
\begin{equation*}
\mid\sqrt{a(\cdot, x,\cdot)}-\sqrt{a(\cdot, y,\cdot)}\mid \leq c\mid x-y\mid^{2+\delta_2},\eqno(2.8)
\end{equation*}
\begin{equation*}
a_{x_i}(\cdot, x,\cdot)=0, x\in \partial\Omega, i=1,2,\cdots, N,\eqno(2.9)
\end{equation*}
\begin{equation*}
  f_i(x)=0, x\in\partial\Omega, i=1,2,\cdots, N.\eqno(2.10)
\end{equation*}
If  $u(x,t)$ and $v(x,t)$ are two solutions of equation (1.1) with the
different initial values $u_{0}(x)$, $v_{0}(x)\in
L^{\infty}(\Omega)$ respectively, then
\begin{equation*}
\int_{\Omega}\mid u(x,t)-v(x,t)\mid dx\leq c\int_{\Omega}\mid
u_{0}(x)-v_{0}(x)\mid dx.\eqno(2.11)
\end{equation*}

Here $d(x)=\text{dist}(x,\partial \Omega)$ is the distance function from the boundary, $a(\cdot, x,t)$ is regarded as the function of the variables $(x,t)$, $a(\cdot, x,\cdot)$ is regarded as the function of $x$.

In general, the conditions listed in Theorem 2.10 are only the sufficient conditions,  and can be replaced by the other assumptions.

If without the condition (2.7), we have

\textbf{Theorem 2.11} Suppose that $a(\cdot,x,t)$ is a $C^{1}(\overline{Q_T})$ function with that $a(\cdot,x,t)\mid_{x\in\partial \Omega}=0$, $a(s,x,t)$ is bounded when $s$ is bounded,
 $f_{i}(x) \in C^{1}(\overline{\Omega})$,  $c(x,t)$ and $g(x,t)$ are bounded. Suppose that the conditions (2.8)-(2.10) are true. If  $u(x,t)$ and $v(x,t)$ are two solutions of equation (1.1) with the different initial values $u_{0}(x)$, $v_{0}(x)\in
L^{\infty}(\Omega)$ respectively, then the stability (2.11) is true.

If the condition (2.10) is not true, we have the following stability based on the partial boundary value condition (1.5) with $\Sigma_p$ appearing as (2.5).

\textbf{Theorem 2.12} Suppose $a(\cdot,x,t)$ is a $C^{1}(\overline{Q_T})$ function, $a(s,x,t)$ is bounded when $s$ is bounded,  $f_{i}(x) \in C^{1}(\overline{\Omega})$,  $c(x,t)$ and $g(x,t)$ are bounded. Suppose that the condition (2.7) is true.
If  $u(x,t)$ and $v(x,t)$ are two solutions of equation (1.1) with the
different initial values $u_{0}(x)$, $v_{0}(x)\in
L^{\infty}(\Omega)$ respectively, and with the same partial boundary value condition
$$
u(x,t)=v(x,t)=0, (x,t)\in \Sigma_p\times (0,T),\eqno(2.12)
$$
then the stability (2.11) is true. Where $\Sigma_p$ has the form (2.5).

Now, we give a simple comment on Theorem 2.11 and  Theorem 2.12. For the linear degenerate parabolic equation
$$
\frac{\partial u}{\partial t} -\frac{\partial }{\partial x_i}\left( a(x)\frac{\partial u}{\partial x_i}\right)
-f_i(x)D_iu+c(x,t)u=g(x,t),\ \ (x,t)\in \Omega \times (0,T),\eqno(2.13)
$$
let $\vec{n}=\{n_i\}$ be the inner normal vector of $\Omega$. To ensure the well-posedness of the solutions of equation (2.13), only a partial boundary value condition (1.5) should be imposed,  where the part of the boundary $\Sigma_p$ can be expressed by Fichera function [26-27]
$$
\Sigma_p=\{x\in \partial \Omega: a(x)>0\}\bigcup\{x\in \partial \Omega: f_i(x)n_i(x)<0\}. \eqno(2.14)
$$

If $a(x)\mid_{x\in\partial \Omega}=0$, and the condition (2.10) is imposed, by (2.14), we have
$$
\Sigma_p=\emptyset.
$$
In the other words, the stability of the weak solutions of equation (2.13) can be obtained independent of the boundary value condition. This coincides with Theorem 2.11.

 If without the condition (2.10), since $a(x)\mid_{x\in\partial \Omega}=0$, (2.14) reduces to the expression (2.5). This coincides with Theorem 2.12.

Fourthly, we would like to suggest that there are many domains satisfying the condition (2.7). For examples,

i) The $N-$dimensinoal cube
$$
C_1=\{x\in\mathbb{R}^N: 0<x_i<1, i=1, 2, \cdots, N\}
$$
the distance function $d$ from the boundary satisfies that
when $x$ is near to the hyperplane $\{x: x_i=0\}$,
$$
d(x)=x_i,
$$
while  $x$ is near to the hyperplane $\{x: x_i=1\}$,
$$
d(x)=1-x_i.
$$

ii) The $N-$ dimensional unit disc
$$
D_1=\{x\in\mathbb{R}^N: |x|<1\}
$$
the distance function  from the boundary is
$$
d(x)=1-r, \ r^2=x_1^2+x_2^2+\cdots+x_N^2,
$$
$$
d_{x_i}=-\frac{x_i}{r},
$$
$$
\Delta d=-\frac{N-1}{r}<0.
$$

The last but not the least, we have said before the condition (2.7) is not a necessary condition. For example, in Theorem 2.11, we have used the condition $a(\cdot,x,t)|_{x\in \partial \Omega}=0$ to replace the condition (2.7). This is very interesting phenomena. Condition (2.7), $\Delta d<0$ reflects the geometric characteristic of the domain $\Omega$, while, $a(\cdot,x,t)$ itself is the diffusion coefficient, the condition $a(\cdot,x,t)|_{x\in \partial \Omega}=0$ implies the diffusion process ends at the boundary $\partial\Omega$. The results of our paper show that these two different conditions both are enough to make the solutions stable.

\bigskip
\section{The proof of the existence}
\bigskip

The existence of the entropy solutions of equation (1.1) can be proved by the similar way as that in [18, 19, 23], we only give the outline of the proof in what follows.

\textbf{Lemma 3.1} [24]  Assume that $\Omega \subset
\mathbb{R}^{N}$ is an open bounded set and  $f_{k},f\in
L^{q}(\Omega)$, as $k\rightarrow \infty $, $f_{k}\rightharpoonup f$
weakly in $L^{q}(\Omega)\ (1\leq q<\infty)$. Then we have
$$
\underset{k\rightarrow \infty}{\lim}\inf \;\parallel f_{k}\parallel
_{L^{q}(\Omega)}^{q}\geq \parallel f\parallel_{L^{q}(\Omega
)}^{q}.
$$

\textbf{Proof of Theorem 2.8}
Consider the regularized problem
\begin{equation*}
\frac{\partial u}{\partial t}=\frac{\partial}{\partial x_i}\left(a(u,x,t)\frac{\partial u}{\partial x_i}\right)+\varepsilon \Delta
u+f_i(x)D_iu-c(x,t)u+g(x,t),\ \ (x,t)\in Q_{T},\eqno(3.1)
\end{equation*}
with the initial-boundary conditions
$$
u(x,0)= u_{0\varepsilon}(x),\ \ x\in \Omega,\eqno(3.2)
$$
$$
u(x,t)= 0,\ \ (x,t)\in \partial \Omega\times
(0,T). \eqno(3.3)
$$
Here, $u_{0\varepsilon}(x)$ is a mollified function of $u_0$. We know
that there exists a classical solutions $u_{\varepsilon}$,  provided
that both $a(u,x,t)$ and $b_{i}(u,x,t)$ satisfy the assumptions given in Theorem 2.8. For more details, one can
refer to [5] or Chapter 8 of [25]. Moreover, we  have
\begin{equation*}
\mid u_{\varepsilon }\mid \leq \parallel u_{0}\parallel_{L^{\infty
}}\leq c.\eqno(3.4)
\end{equation*}

\textbf{Step 1} Multiplying equation
(2.1) with $u_{\varepsilon}$,  it is easy to show that
\begin{equation*}
\iint_{Q_T}a(u_{\varepsilon}, x,t)|\nabla u_{\varepsilon}|^2dxdt\leq c.\eqno(3.5)
\end{equation*}

Then, $\sqrt{a(u_{\varepsilon},x,t)}\frac{\partial u_{\varepsilon}}{\partial x_{i}}$
is weakly compact in $L^2(Q_T)$. By choosing a subsequence (still denoting it as $\sqrt{a(u_{\varepsilon},x,t)}\frac{\partial u_{\varepsilon}}{\partial x_{i}}$) , we are able to show that
$$
\sqrt{a(u_{\varepsilon},x,t)}\frac{\partial u_{\varepsilon}}{\partial x_{i}}\rightharpoonup \widehat{\sqrt{a(u,x,t)}}\frac{\partial u}{\partial x_i}, \text{in}\  L^2(Q_T),
$$
$u$ satisfies (1) of Definition 2.6.

\textbf{Step 2} Let $\varphi\in C_0^{2}({Q_{T}}),\;\varphi\geq
0$. Multiplying both sides of (2.1) by
$\varphi S_{\eta}(u_{\varepsilon}-k)$,
  integrating it by part, we can deduce that
$$
\begin{gathered}
\iint\nolimits_{Q_{T}}I_{\eta}(u_{\varepsilon}-k)\varphi
_{t}dxdt+\iint\nolimits_{Q_{T}}A_{\eta}(u_{\varepsilon},x,t,k)\triangle\varphi
dxdt
\hfill \\
 \begin{array}{*{20}{c}}
{}
\end{array}
-\iint\nolimits_{Q_{T}}I_{\eta}(u_{\varepsilon}-k)f_{i}(x)\varphi
_{x_{i}}dxdt
-\varepsilon\iint\nolimits_{Q_{T}}\nabla
u_{\varepsilon}\cdot\nabla\varphi S_{\eta}(u_{\varepsilon
}-k)dxdt
\hfill \\
 \begin{array}{*{20}{c}}
{}
\end{array}
-\varepsilon\iint\nolimits_{Q_{T}}\mid\nabla
u_{\varepsilon}\mid^{2}S_{\eta}^{\prime}(u_{\varepsilon
}-k)\varphi dxdt
+\iint\nolimits_{Q_{T}}\int_{k}^{u_{\varepsilon}}a_{x_i}(s,x,t)S_{\eta}(s-k)ds\varphi_{x_{i}}dxdt
\hfill \\
 \begin{array}{*{20}{c}}
{}
\end{array}
 -\iint\nolimits_{Q_{T}}a(u_{\varepsilon},x,t)\mid\nabla u_{\varepsilon}\mid^{2}S_{\eta}^{\prime}(u_{\varepsilon}-k)\varphi dxdt
\hfill \\
 \begin{array}{*{20}{c}}
{}
\end{array}
-\iint\nolimits_{Q_{T}}f_{ix_i}(x)(u_{\varepsilon}-k)\varphi S_{\eta
}(u_{\varepsilon}-k)dxdt+\iint\nolimits_{Q_{T}}f_{ix_i}(x)\int_{k}^{u_{\varepsilon}}(s-k)h_{\eta}(s-k)ds\varphi dxdt
\hfill \\
 \begin{array}{*{20}{c}}
{}
\end{array}
-\iint \nolimits_{Q_{T}}[c(x,t)u_{\varepsilon}+g(x,t)]\varphi S_{\eta}(u_{\varepsilon}-k)dxdt
\hfill \\
 \begin{array}{*{20}{c}}
{}
\end{array}
= 0.
\hfill\\
\end{gathered}
\eqno(3.6)
$$
By Lemma 3.1, we have
$$
 \underset{\varepsilon \rightarrow 0}{\lim\inf
}\iint\nolimits_{Q_{T}}S_{\eta}^{\prime }(u_{\varepsilon
}-k)a(u_{\varepsilon},x,t)\frac{\partial u_{\varepsilon}}{\partial
x_{i}}\frac{\partial u_{\varepsilon}}{\partial x_{i}}\varphi
dxdt
\geq  \sum_{i=1}^{N}\iint\nolimits_{Q_{T}}S_{\eta}^{\prime}(u-k)\mid g^{i}\mid^{2}\varphi dxdt.\eqno(3.7)
$$
Letting $\varepsilon \rightarrow 0$ in (3.6), it is easily to obtain (2.3).

\textbf{Step 3} At last,  the initial value (1.2) is true in the sense of (2.4), its proof can be found in [21].

Thus, the existence of the entropy solution in the sense of Definition 2.6 has been proved,  Theorem 2.8 follows immediately.

\textbf{Lemma 3.2} Let $u_{\varepsilon}$ be the solution of the problem (3.1)-(3.3).
If the assumptions given in Theorem 2.9 hold,
then
\begin{equation*}
|\text{grad}u_{\varepsilon}|_{L^{1}(\Omega)}\leq c,
\end{equation*}
where $c$ is independent of $\varepsilon$, and
$$ |\text{grad}u_{\varepsilon}|^{2}=\sum_{i=1}^{N} \left| \frac{\partial u_{\varepsilon}}{\partial x_{i}} \right|^{2}
+\left|\frac{\partial u}{\partial t} \right|^{2}.$$

Lemma  3.2 can be proved in a similar manner as Theorem 11 of [29], we omit the details here.

By Theorem 2.8 and Lemma 3.2, we know that Theorem 2.9 is true.

\bigskip
\section{ Kruzkov's bi-variables method}
\bigskip

Similar as [1, 22], we denote that $\Gamma_{u}$ is the set of all jump points of $u\in BV(Q_{T})$, $v$ is
the normal of $\Gamma_{u}$ at $X=(x,t)$,  $u^{+}(X)$ and $u^{-}(X)$ are
the approximate limits of $u$ at $X\in \Gamma_{u}$ with respect to
$(v,Y-X)>0$ and $(v,Y-X)<0$, respectively. For the continuous functions
$p(u,x,t)$ and $u\in BV(Q_{T})$, the composite mean value of $p$ is defined as
$$
\widehat{p}(u,x,t)=\int_{0}^{1}p(\tau u^{+}+(1-\tau )u^{-},x,t)d\tau.
$$
If
$f(s)\in C^{1}( \mathbb{R})$ and $u\in BV(Q_{T})$, then $f(u)\in
BV(Q_{T})$ and
$$
\frac{\partial f(u)}{\partial x_{i}}=\widehat{f^{\prime
}}(u)\frac{\partial u}{\partial x_{i}},\ \ \;i=1, 2, \cdots, N, N+1,
$$
where $x_{N+1}=t$.

\textbf{Lemma 4.1} Let $u$ be a solution of (1.1). Then
\begin{equation*}
a(s,x,t)=0,\;s\in I(u^{+}(x,t),u^{-}(x,t))\;\;a.e.\;on\;\Gamma
_{u},\eqno(4.1)
\end{equation*}
which $I(\alpha,\beta)$ denote the closed interval with endpoints
$\alpha $ and $\beta$, and (4.1) is in the sense of Hausdorff
measure $H_{N}(\Gamma_{u})$.

\textbf{Proof}\ Denote
\begin{equation*}
\Gamma _{1}=\{(x,t)\in \Gamma_{u}, v_{1}(x,t)=\cdots =v_{N}(x,t)=0\},
\end{equation*}
\begin{equation*}
\Gamma_{2}=\{(x,t)\in \Gamma_{u}, v_{1}^{2}(x,t)+\cdots
+v_{N}^{2}(x,t)>0\}.
\end{equation*}

At first, we prove $a(s,x,t)=0,\;s\in
I(u^{+}(x,t),u^{-}(x,t))\;\;a.e.\;on\;\Gamma_{1}$. Since any
measurable subset of $\Gamma_{1}$ can be expressed as the union of
Borel sets and a set of measure zero, it suffices to prove
\begin{equation*}
a(s)=0,\;s\in I(u^{+}(x,t),u^{-}(x,t))\;\;a.e.\;on\;U\subset \Gamma
_{1},
\end{equation*}
where $U$ is a Borel subset of $\Gamma_{1}$.  For any bounded
function $f(x,t)$, which is measurable with respect to measure
$\frac{\partial u}{
\partial x_{i}}$,  Lemma 3.7.8 in [1] shows that
\begin{equation*}
\iint\nolimits_{U}f(x,t)\frac{\partial u}{\partial x_{i}}
=\int_{0}^{T}dt\int_{U^{t}}f(x,t)\frac{\partial u}{\partial
x_{i}},\eqno(4.2)
\end{equation*}
\ where $U^{t}=\{x:(x,t)\in U\}$. Moreover, for any Borel subset
$S\subset U$, $S^{t}\subset U^{t}$, for $i=1, 2, \cdots, N$,
$$
\frac{\partial u}{\partial x_{i}}(S)
=\int_{S}(u^{+}(x,t)-u^{-}(x,t))v_{i}dH,
$$
$$
\frac{\partial u(\cdot,t)}{\partial x_{i}}(S^{t})
=\int_{S^{t}}(u_{+}^{t}(x,t)-u_{-}^{t}(x,t))v_{i}dH^{t}.
$$
(4.2) is equivalent to
\begin{equation*}
\iint\nolimits_{U}f(x,t)(u^{+}(x,t)-u^{-}(x,t))v_{i}dH=\int_{0}^{T}dt
\int_{U^{t}}f(x,t)(u_{+}^{t}(x,t)-u_{-}^{t}(x,t))v_{i}^{t}dH^{t}.\
\end{equation*}
The definition of $\Gamma_{1}$ implies that the left hand side
vanishes, then
\begin{equation*}
\int_{0}^{T}dt
\int_{U^{t}}f(x,t)(u_{+}^{t}(x,t)-u_{-}^{t}(x,t))v_{i}^{t}dH^{t}=0.
\end{equation*}
If we choose $f(x,t)=\chi_{u}(x,t)$sgn$(u_{+}^{t}(x,t)-u_{-}^{t}(x,t))$
sgn$v_{i}^{t}$, where $\chi _{u}(x,t)$ is the characteristic
function of $U$,  sum up for $i$ from $1$ up to $N$, then
\begin{equation*}
\int_{G}dt\int_{U^{t}}(u_{+}^{t}(x,t)-u_{-}^{t}(x,t))(\mid
v_{1}^{t}\mid +\cdots +\mid v_{N}^{t}\mid)dH^{t}=0,
\end{equation*}
where $G$ is the projection of $U$ on the $t$-axis. (4.2) implies
for almost all $t\in G$,
\begin{equation*}
\int_{U^{t}}(u_{+}^{t}(x,t)-u_{-}^{t}(x,t))(\mid v_{1}^{t}\mid
+\cdots +\mid v_{N}^{t}\mid )dH^{t}=0,
\end{equation*}
and hence for almost all $t\in G$,
\begin{equation*}
v_{1}^{t}=\cdots =v_{N}^{t}=0,
\end{equation*}
$H^{t}$-almost everywhere on $U^{t}$, which is impossible unless
mes$G=0$.

For any $\alpha,\beta $ with $0<\alpha <\beta <T$, we choose $\psi
_{j}(t)\in C_{0}^{\infty}(0,T)$ such that
\begin{equation*}
0\leq \psi_{j}(t)\leq 1,\underset{j\rightarrow \infty}{\lim}\psi
_{j}(t)=\chi_{\lbrack \alpha,\beta]}(t),\;\forall t\in \lbrack 0,T],
\end{equation*}
and choose $\varphi_{n}\in C_{0}^{\infty}(Q_{T})$ such
that
\begin{equation*}
\mid \varphi_{n}(x,t)\mid \leq 1,\;\underset{n\rightarrow \infty
}{\lim}\varphi_{n}=\chi_{U}\;\text{in}\ L^{1}(Q_{T},\mid
\frac{\partial u}{\partial t}\mid).
\end{equation*}
Now, denoting that
$$
A(u,x,t)=\int_0^ua(s,x,t)ds,
$$
 from the definition of BV-function, we have
\begin{eqnarray*}
&&\iint\nolimits_{Q_{T}}\varphi_{n}(x,t)\psi_{j}(t)\frac{\partial
u}{
\partial t}
\\
&&=\iint\nolimits_{Q_{T}}A(u,x,t)\Delta \varphi_{n}(x,t)\psi
_{j}(t)dxdt-\iint\nolimits_{Q_{T}}a_{x_i}(s,x,t)ds\varphi_{nx_i}(x,t)\psi_{j}(t)dxdt
\\
&&-\iint\nolimits_{Q_{T}}f_{i}(x)u\frac{\partial}{\partial
x_{i}}\varphi_{n}(x,t)\psi_{j}(t)dxdt+\iint\nolimits_{Q_{T}}f_{ix_i}u\varphi_{n}(x,t)\psi_{j}(t)
\\
&&+\iint\nolimits_{Q_{T}}[g(x,t)-c(x,t)u]\varphi_{n}(x,t)\psi_{j}(t)dxdt.
\end{eqnarray*}
Let $j\rightarrow \infty $. Then
\begin{eqnarray*}
&&\iint\nolimits_{Q_{T}}\varphi_{n}(x,t)\chi_{\lbrack \alpha,\beta
]}(t)\frac{\partial
u}{
\partial t}
\\
&&=\iint\nolimits_{Q_{T}}A(u,x,t)\Delta \varphi_{n}(x,t)\chi_{\lbrack \alpha,\beta
]}(t)dxdt-\iint\nolimits_{Q_{T}}a_{x_i}(s,x,t)ds\varphi_{nx_i}(x,t)\chi_{\lbrack \alpha,\beta
]}(t)dxdt
\\
&&-\iint\nolimits_{Q_{T}}f_{i}(x)u\frac{\partial}{\partial
x_{i}}\varphi_{n}(x,t)\chi_{\lbrack \alpha,\beta
]}(t)dxdt+\iint\nolimits_{Q_{T}}f_{ix_i}u\varphi_{n}(x,t)\chi_{\lbrack \alpha,\beta
]}(t)
\\
&&+\iint\nolimits_{Q_{T}}[g(x,t)-c(x,t)u]\varphi_{n}(x,t)\chi_{\lbrack \alpha,\beta
]}(t)dxdt.
\end{eqnarray*}
Clearly, this equality also holds if $[\alpha,\beta]$ is replaced by
$(\alpha,\beta )$ and hence it holds even if $[\alpha,\beta]$ is
replaced by any open set $I$ with $\overline{I}\subset (0,T)$. Since
$G$ is a Borel set, by approximation we may conclude that
\begin{eqnarray*}
&&\iint\nolimits_{Q_{T}}\varphi_{n}(x,t)\chi_{G}(t)\frac{\partial
u}{
\partial t}
\\
&&=\iint\nolimits_{Q_{T}}A(u,x,t)\Delta \varphi_{n}(x,t)\chi_{G}(t)dxdt-\iint\nolimits_{Q_{T}}a_{x_i}(s,x,t)ds\varphi_{nx_i}(x,t)\chi_{G}(t)dxdt
\\
&&-\iint\nolimits_{Q_{T}}f_{i}(x)u\frac{\partial}{\partial
x_{i}}\varphi_{n}(x,t)\chi_{G}(t)dxdt+\iint\nolimits_{Q_{T}}f_{ix_i}u\varphi_{n}(x,t)\chi_{G}(t)
\\
&&+\iint\nolimits_{Q_{T}}[g(x,t)-c(x,t)u]\varphi_{n}(x,t)\chi_{G}(t)dxdt.
\end{eqnarray*}

The two terms on the right hand vanish by that mes$G=0$, and
\begin{equation*}
\iint\nolimits_{Q_{T}}\varphi_{n}(x,t)\chi_{G}(t)\frac{\partial
u}{\partial t}=0.
\end{equation*}
Let $n\rightarrow \infty $. Then
\begin{equation*}
\iint\nolimits_{U}\frac{\partial u}{\partial
t}=\iint\nolimits_{Q_{T}}\chi_{U}(x,t)\chi_{G}\frac{\partial
u}{\partial t}=0.
\end{equation*}
Hence
\begin{equation*}
\int_{U}(u^{+}(x,t)-u^{-}(x,t))v_{t}dH=0,
\end{equation*}
which implies $H(U)=0$ and $H(\Gamma_{1})=0$ by the arbitrariness of
$U$.

Secondly, we prove $H(\Gamma_{2})=0$. Let $U$ be any Borel subset of
$\Gamma_{2}$ which is compact in $Q_{T}$. Since $U$ is a set of
$N+1$-dimensional measure zero and $\frac{\partial}{\partial x_i} A(u,x,t)\in
L_{loc}^{2}(Q_{T})$, we have
\begin{equation*}
\iint\nolimits_{U}\frac{\partial}{\partial
x_{i}}A(u,x,t)dxdt=0,\;i=1,\cdots, N,
\end{equation*}
and hence
\begin{equation*}
\int_{U}[A(u^{+},x,t)-A(u^{-}, x,t)]v_{i}dH=0,\;i=1,\cdots, N.
\end{equation*}
Form this fact, it follows by the definition of $\Gamma_{2}$ that
\begin{equation*}
\int_{u^{-}(x,t)}^{u^{+}(x,t)}a(s,x,t)ds=0,\;a.e.\ on\ \Gamma_{2}.
\end{equation*}
Thus the lemma is proved.

In this section, we apply Kru\v{z}kov bi-variables method to the main equation (1.1). In details, let $u(x,t)$ and $v(x,t)$ be two entropy solutions
of equation (1.1) with the initial values
\begin{equation*}
u(x,0)=u_{0}(x)\ \ \text{and}\ \ v(x,0)=v_{0}(x),
\end{equation*}
respectively.

By Definition 2.6, for any nonnegative $\varphi\in C_{0}^{2}(Q_{T})$,   we have
$$
\begin{gathered}
\iint\nolimits_{Q_{T}}\left[I_{\eta}(u-k)\varphi_{t}-f_i(x)I_{\eta}(u-k)\varphi_{x_{i}}+A_{\eta}(u,x,t,k)\Delta \varphi
-\sum_{i=1}^{N}S_{\eta}^{\prime}(u-k)\mid g^{i}\mid
^{2}\varphi\right]dxdt
\hfill \\
 \begin{array}{*{20}{c}}
{}
\end{array}
+\iint \nolimits_{Q_{T}}\int_{k}^{u}a_{x_i}(s,x,t)S_{\eta}(s-k)ds\varphi_{x_{i}}dxdt
\hfill \\
 \begin{array}{*{20}{c}}
{}
\end{array}
-\iint\nolimits_{Q_{T}}f_{ix_i}(x)(u-k)\varphi S_{\eta
}(u-k)dxdt+\iint\nolimits_{Q_{T}}f_{ix_i}(x)\int_{k}^{u}(s-k)h_{\eta}(s-k)ds\varphi dxdt
\hfill \\
 \begin{array}{*{20}{c}}
{}
\end{array}
-\iint \nolimits_{Q_{T}}[c(x,t)u-g(x,t)]\varphi S_{\eta}(u-k) dxdt
\hfill \\
 \begin{array}{*{20}{c}}
{}
\end{array}
\geq 0,
\hfill\\
\end{gathered}
\eqno(4.3)
$$
and
$$
\begin{gathered}
\iint\nolimits_{Q_{T}}\left[I_{\eta}(v-l)\varphi_{\tau}-f_i(y)I_{\eta}(v-l)\varphi_{y_{i}}+A_{\eta}(v,y.\tau,l)\Delta \varphi
-\sum_{i=1}^{N}S_{\eta}^{\prime}(v-l)\mid g^{i}\mid
^{2}\varphi\right]dyd\tau
\hfill \\
 \begin{array}{*{20}{c}}
{}
\end{array}
+\iint \nolimits_{Q_{T}}\int_{l}^{v}a_{y_i}(s,y,\tau)S_{\eta}(v-l)ds\varphi_{y_{i}}dyd\tau
\hfill \\
 \begin{array}{*{20}{c}}
{}
\end{array}
-\iint\nolimits_{Q_{T}}f_{iy_i}(y)(v-l)\varphi S_{\eta
}(v-l)dyd\tau+\iint\nolimits_{Q_{T}}f_{iy_i}(y)\int_{l}^{v}(s-k)h_{\eta}(s-l)ds\varphi dyd\tau
\hfill \\
 \begin{array}{*{20}{c}}
{}
\end{array}
-\iint \nolimits_{Q_{T}}[c(y,\tau)v-g(y,\tau)]\varphi S_{\eta}(v-l) dyd\tau
\geq 0.
\hfill\\
\end{gathered}
\eqno(4.4)
$$
Let
$$
\psi (x,t,y,\tau)=\phi (x,t)j_{h}(x-y,t-\tau ),
$$
for any $\phi (x,t)\geq 0,\;\phi (x,t)\in C_{0}^{\infty}(Q_{T})$, and
\begin{equation*}
j_{h}(x-y,t-\tau )=\omega_{h}(t-\tau)\Pi_{i=1}^{N}\omega
_{h}(x_{i}-y_{i}).
\end{equation*}
Here,
$ \omega _{h}(s)=\frac{1}{h}\omega(\frac{s}{h}),\ \omega(s)\in
C_{0}^{\infty}(\mathbb{R}),\ \omega(s)\geq 0,\ \omega (s)=0\ \text{if}\ \mid s\mid
>1,$ and $\int_{-\infty}^{\infty}\omega (s)ds=1.$
Moreover, for any given positive constant $\delta$, there holds
\begin{equation*}
\lim_{h\rightarrow 0}\omega _{h}^{\prime }(s)s^{2+\delta }=0.\eqno(4.5)
\end{equation*}

We choose $k=v(y,\tau),\;l=u(x,t)$ and $\varphi=\psi (x,t,y,\tau)$
in (4.3) and (4.4). Integrating it over $Q_{T}$, using the fact of that $S_{\eta}(u-v)=-S_{\eta}(v-u)$, we have
$$
\begin{gathered}
 \iint\nolimits_{Q_{T}}\iint\nolimits_{Q_{T}}\bigg\{I_{\eta
}(u-v)(\psi _{t}+\psi_{\tau})+A_{\eta}(u,x,t,v)\Delta_{x}\psi+A_{\eta
}(v,y,\tau,u)\Delta_{y}\psi
\hfill \\
 \begin{array}{*{20}{c}}
{}
\end{array}
+\int_{v}^{u}a_{x_i}(s,x,t)S_{\eta}(s-v)ds\psi_{x_{i}}
+\int_{u}^{v}a_{y_i}(s,y,\tau)S_{\eta}(s-u)ds\psi_{y_{i}}
\hfill \\
 \begin{array}{*{20}{c}}
{}
\end{array}
-\sum_{i=1}^{N}S_{\eta}^{\prime}(u-v)\left[\mid g^{i}(u,x,t)\mid ^{2}+\mid g^{i}(v,y,\tau)\mid ^{2}\right]\psi
\hfill \\
 \begin{array}{*{20}{c}}
{}
\end{array}
-\left[f_i(x)I_{\eta}(u-v)\psi_{x_{i}}+f_i(y)I_{\eta}(v-u)\psi_{y_{i}}\right]
\hfill \\
 \begin{array}{*{20}{c}}
{}
\end{array}
-\left[\text{div}_{x}\vec{f}S_{\eta}(u-v)(u-v)+\text{div}_{y}\vec{f}S_{\eta}(v-u)(v-u)\right]\psi
\hfill \\
 \begin{array}{*{20}{c}}
{}
\end{array}
+\left[\text{div}_{x}\vec{f}\int_{v}^{u}(s-v)h_{\eta}(s-v)ds+\text{div}_{y}\vec{f}\int_{u}^{v}(s-u)h_{\eta}(s-u)ds\right]\psi
\hfill \\
 \begin{array}{*{20}{c}}
{}
\end{array}
-[c(x,t)u-c(y,\tau)v]S_{\eta}(u-v)\psi+[(g(x,t)-g(y,\tau)]S_{\eta}(u-v)\psi \bigg\}dxdtdyd\tau
\hfill \\
 \begin{array}{*{20}{c}}
{}
\end{array}
\geq 0.
\hfill\\
\end{gathered}
\eqno(4.6)
$$
We can use the facts
\begin{eqnarray*}
&&  \frac{\partial j_{h}}{\partial t}+\frac{\partial j_{h}}{\partial
\tau }=0,\ \ \frac{\partial j_{h}}{\partial x_{i}}+\frac{\partial
j_{h}}{\partial y_{i}}=0,\ \ i=1,\cdots ,N,\\
&& \frac{\partial \psi}{\partial t}+\frac{\partial \psi}{\partial \tau
}=\frac{\partial \phi}{\partial t}j_{h},\ \ \frac{\partial \psi
}{\partial x_{i}}+\frac{\partial \psi}{\partial
y_{i}}=\frac{\partial \phi}{\partial x_{i}}j_{h},
\end{eqnarray*}
to analysis every term of the left hand side of (4.6).

The first term, we have
$$
\lim_{h \rightarrow 0}\lim_{\eta \rightarrow 0}\iint\nolimits_{Q_{T}}\iint\nolimits_{Q_{T}}I_{\eta
}(u-v)\psi _{t}dxdtdyd\tau=\iint\nolimits_{Q_{T}}|u(x.t)-v(x,t)|\phi_{t}dxdt.\eqno(4.7)
$$
From the second term to the sixth term, by a very complicated calculations [23], by (4.1) in Lemma 4.1, using the condition (2.8) and the observation (4.5), we can deduce that
$$
\begin{gathered}
\lim_{h \rightarrow 0}\lim_{\eta \rightarrow 0}\iint\nolimits_{Q_{T}}\iint\nolimits_{Q_{T}}\bigg\{I_{\eta
}(u-v)(\psi _{t}+\psi_{\tau})+A_{\eta}(u,x,t,v)\Delta_{x}\psi+A_{\eta
}(v,y,\tau,u)\Delta_{y}\psi
\hfill \\
 \begin{array}{*{20}{c}}
{}
\end{array}
+\int_{v}^{u}a_{x_i}(s,x,t)S_{\eta}(s-v)ds\psi_{x_{i}}
+\int_{u}^{v}a_{y_i}(s,y,\tau)S_{\eta}(s-u)ds\psi_{y_{i}}
\hfill \\
 \begin{array}{*{20}{c}}
{}
\end{array}
-\sum_{i=1}^{N}S_{\eta}^{\prime}(u-v)\left[\mid g^{i}(u,x,t)\mid ^{2}+\mid g^{i}(v,y,\tau)\mid ^{2}\right]\psi\bigg\}dxdtdyd\tau
\hfill \\
 \begin{array}{*{20}{c}}
{}
\end{array}
=\iint\nolimits_{Q_{T}}\bigg\{\text{sgn}(u-v)(A(u,x,t)-A(v,x,t))\Delta\phi
\hfill \\
 \begin{array}{*{20}{c}}
{}
\end{array}
+\int_{v}^{u}a_{x_{i}}(s,x,t)\text{sgn}(s-v)ds\phi _{x_{i}}+\int_{u}^{v}a_{x_{i}}(s,x,t)\text{sgn}(s-u)ds\phi
_{x_{i}}\bigg\}dxdt.
\hfill\\
\end{gathered}
\eqno(4.8)
$$

For the seventh term, by the fact
$$
\psi_{y_i}=\phi_{x_i}j_h-\psi_{x_i},
$$
we have
$$
\begin{gathered}
\lim_{h \rightarrow 0}\lim_{\eta\rightarrow 0}\iint\nolimits_{Q_{T}}\iint\nolimits_{Q_{T}}\left[f_i(x)I_{\eta}(u-k)\psi_{x_{i}}+f_i(y)I_{\eta}(v-l)\psi_{y_{i}}\right]dxdtdyd\tau
\hfill \\
 \begin{array}{*{20}{c}}
{}
\end{array}
=\lim_{h \rightarrow 0}\iint\nolimits_{Q_{T}}\iint\nolimits_{Q_{T}}\left[f_i(x)\psi_{x_{i}}+f_i(y)\psi_{y_{i}}\right]|u-v|dxdtdyd\tau
\hfill \\
 \begin{array}{*{20}{c}}
{}
\end{array}
=\lim_{h \rightarrow 0}\iint\nolimits_{Q_{T}}\iint\nolimits_{Q_{T}}\left[f_i(x)\psi_{x_{i}}+f_i(y)(\phi_{x_i}j_h-\psi_{x_i})\right]|u-v|dxdtdyd\tau
\hfill \\
 \begin{array}{*{20}{c}}
{}
\end{array}
=\iint\nolimits_{Q_{T}}f_i(x)\phi_{x_i}|u-v|dxdt.
\hfill\\
\end{gathered}
\eqno(4.9)
$$

For the eighth term, it is obviously
$$
\begin{gathered}
-\lim_{h \rightarrow 0}\lim_{\eta \rightarrow 0}\iint\nolimits_{Q_{T}}\iint\nolimits_{Q_{T}}\left[\text{div}_{x}\vec{f}S_{\eta}(u-v)(u-v)+\text{div}_{y}\vec{f}S_{\eta}(v-u)(v-u)\right]\psi dxdtdyd\tau
\hfill \\
 \begin{array}{*{20}{c}}
{}
\end{array}
=-\iint\nolimits_{Q_{T}}\text{div}\vec{f}|u-v|\phi dxdt.
\hfill\\
\end{gathered}
\eqno(4.10)
$$

For the ninth term, it is obviously
$$
\begin{gathered}
-\lim_{h \rightarrow 0}\lim_{\eta \rightarrow 0}\iint\nolimits_{Q_{T}}\iint\nolimits_{Q_{T}}\left[\text{div}_{x}\vec{f}\int_{v}^{u}(s-v)h_{\eta}(s-v)ds+\text{div}_{y}\vec{f}\int_u^v(s-u)h_{\eta}(s-u)ds\right]\psi dxdtdyd\tau
\hfill \\
 \begin{array}{*{20}{c}}
{}
\end{array}
=0.
\hfill\\
\end{gathered}
\eqno(4.11)
$$
For the tenth term,
$$
\begin{gathered}
-\lim_{h \rightarrow 0}\lim_{\eta \rightarrow 0}\iint\nolimits_{Q_{T}}\iint\nolimits_{Q_{T}}\left[c(x,t)u-c(y,\tau)v\right]S_{\eta}(u-v)\psi dxdtdyd\tau
\hfill \\
 \begin{array}{*{20}{c}}
{}
\end{array}
=-\iint\nolimits_{Q_{T}}c(x,t)|u-v|\phi dxdt.
\hfill\\
\end{gathered}
\eqno(4.12)
$$
For the last term,
$$
\begin{gathered}
\lim_{h \rightarrow 0}\lim_{\eta \rightarrow 0}\iint\nolimits_{Q_{T}}\iint\nolimits_{Q_{T}}\left[g(x,t)-g(y,\tau)\right]S_{\eta}(u-v)\psi dxdtdyd\tau
\hfill \\
 \begin{array}{*{20}{c}}
{}
\end{array}
=\iint\nolimits_{Q_{T}}[g(x,t)g(x,t)]\text{sign}(u-v)\phi dxdt.
\hfill\\\begin{array}{*{20}{c}}
{}
\end{array}
=0.
\hfill\\
\end{gathered}
\eqno(4.13)
$$
Thus, if we let $\eta \rightarrow 0$ and $h\rightarrow 0$  in (4.6), then we have
$$
\begin{gathered}
 \iint_{Q_{T}}\bigg\{\mid u(x,t)-v(x,t)\mid \phi
_{t}+\text{sgn}(u-v)[A(u,x,t)-A(v,x,t)]\Delta\phi
\hfill \\
 \begin{array}{*{20}{c}}
{}
\end{array}
+\int_{v}^{u}a_{x_{i}}(s,x,t)\text{sgn}(s-v)ds\phi _{x_{i}}+\int_{u}^{v}a_{x_{i}}(s,x,t)\text{sgn}(s-u)ds\phi_{x_{i}}
\hfill \\
 \begin{array}{*{20}{c}}
{}
\end{array}
-[f_i(x)\phi_{x_i}+\text{div}\vec{f}\phi+c(x,t)\phi]|u-v|\bigg\}dxdt
\hfill \\
 \begin{array}{*{20}{c}}
{}
\end{array}
\geq 0.
\hfill\\
\end{gathered}
\eqno(4.14)
$$

By choosing some special test functions or some special domains $\Omega$, one can prove the stability of the entropy solutions according to (4.14).

\bigskip
\section{ The proof of Theorem 2.10 and Theorem 2.11}
\bigskip

\textbf{Proof of Theorem 2.10}
For small enough $\lambda$, we define
$$
\varphi_{\lambda}(x)=\left\{\begin{array}{cc}
-\frac{(d-\lambda)^2}{\lambda^2}+1,\ \ &\ \ {\rm if}\ \ 0\leq d\leq \lambda,\\
1,\ \ &\ \ {\rm if}\ d\geq \lambda.
\end{array}
\right.\eqno(5.1)
$$
By a process of limit, we can choose the test function  in (4.13) as
$$
\phi(x,t)=\eta(t)\varphi_{\lambda}(x),\eqno(5.2)
$$
where $0\leq\eta(t)\in C_0^{1}(t)$.

When $x\in \Omega_{\lambda}=\{x\in \Omega: d(x)<\lambda\}$,
$$
\partial_{x_{i}}\phi(x,t)=\eta(t)\partial_{x_{i}}\varphi_{\lambda}(x)=-\eta(t)\frac{2(d-\lambda)}{\lambda^2}d_{x_{i}},
$$
$$
\Delta\phi=-\eta(t)\left[\frac{2}{\lambda^2}|\nabla d|^2+\frac{2(d-\lambda)}{\lambda^2}\Delta d\right].
$$
While in $\Omega\setminus\Omega_{\lambda}$,
$$
\phi_{x_i}=0, \Delta \phi=0.
$$

In the first place, by the assumption of that $\Delta d\leq 0$, choosing $\lambda$ is small enough,  when $x$ is near to the boundary, $d(x)<\lambda$, we have
$$
\begin{gathered}
\int_{\Omega}\text{sgn}(u-v)(A(u,x,t)-A(v,x,t))\Delta\phi dx
\hfill \\
 \begin{array}{*{20}{c}}
{}
\end{array}
=-\eta(t)\int_{\Omega_{\lambda}}|A(u,x,t)-A(v,x,t)|\left[\frac{2}{\lambda^2}|\nabla d|^2+\frac{2(d-\lambda)}{\lambda^2}\Delta d\right] dx
\hfill \\
 \begin{array}{*{20}{c}}
{}
\end{array}
\leq 0.
\hfill\\
\end{gathered}
\eqno(5.3)
$$

In the second place, by that $|d_{x_i}|\leq |\nabla d|=1$, and by (2.9),  $a_{x_i}(s,x,t)=0$ when $x\in \partial \Omega$,
$$
\begin{gathered}
\lim_{\lambda\rightarrow 0}\left|\int_{\Omega}\int_{v}^{u}a_{x_{i}}(s,x,t)\text{sgn}(s-v)ds\phi _{x_{i}}dx\right|
\hfill \\
 \begin{array}{*{20}{c}}
{}
\end{array}
=\lim_{\lambda\rightarrow 0}\left|\int_{\Omega_{\lambda}}\int_{v}^{u}a_{x_{i}}(s,x,t)\text{sgn}(s-v)ds\phi _{x_{i}}dx\right|
\hfill \\
 \begin{array}{*{20}{c}}
{}
\end{array}
\leq \lim_{\lambda\rightarrow 0}\frac{c}{\lambda}\int_{\Omega_{\lambda}}|\int_{v}^{u}a_{x_{i}}(s,x,t)\text{sgn}(s-v)ds|dx
\hfill \\
 \begin{array}{*{20}{c}}
{}
\end{array}
=\int_{\partial\Omega}|\int_{v}^{u}a_{x_{i}}(s,x,t)\text{sgn}(s-v)ds|d\Sigma
\hfill \\
 \begin{array}{*{20}{c}}
{}
\end{array}
=0.
\hfill\\
\end{gathered}
\eqno(5.4)
$$

Similarly, we have
$$
\lim_{\lambda\rightarrow 0}\left|\int_{\Omega}\int_{u}^{v}a_{x_{i}}(s,x,t)\text{sgn}(s-u)ds\phi_{x_{i}}dx\right|=0.\eqno(5.5)
$$

Moreover, by that $|d_{x_i}|\leq |\nabla d|=1$, and by the assumption of that $f_i(x)=0$ when $x\in\partial\Omega$, we have
$$
\begin{gathered}
\lim_{\lambda\rightarrow 0}\left|\int_{\Omega}f_{i}(x)\phi_{x_{i}}(u-v)dx\right|
\hfill \\
 \begin{array}{*{20}{c}}
{}
\end{array}
\leq 2\lim_{\lambda\rightarrow 0}\int_{\Omega_{\lambda}}|f_{i}(x)|\frac{|(d-\lambda)d_{x_i}|}{\lambda^2}\eta(t)|u-v|dx
\hfill \\
 \begin{array}{*{20}{c}}
{}
\end{array}
\leq c\sum_{i=1}^N\lim_{\lambda\rightarrow 0}\frac{1}{\lambda}\int_{\Omega_{\lambda}}|f_{i}(x)|\eta(t)|u-v|dx
\hfill \\
 \begin{array}{*{20}{c}}
{}
\end{array}
=c\sum_{i=1}^N\int_{\partial\Omega}|f_{i}(x)|\eta(t)|u-v|d\Sigma
\hfill \\
 \begin{array}{*{20}{c}}
{}
\end{array}
=0,
\hfill\\
\end{gathered}
\eqno(5.6)
$$
and it is clearly that
$$
\begin{gathered}
\lim_{\lambda\rightarrow 0}\int_{\Omega}[\text{div}\vec{f}+c(x,t)]\phi|u-v|dx
\hfill \\
 \begin{array}{*{20}{c}}
{}
\end{array}
=\int_{\Omega}[\text{div}\vec{f}+c(x,t)]\eta(t)|u-v|dx
\hfill \\
 \begin{array}{*{20}{c}}
{}
\end{array}
\leq c\int_{\Omega}\eta(t)|u-v|dx.
\hfill\\
\end{gathered}
\eqno(5.7)
$$

By (5.3)-(5.7), according to (4.12), we have
$$
\iint_{Q_{T}}|u(x,t)-v(x,t)|\phi_{t}dxdt
+c\int_{0}^{T}\int_{\Omega}\eta(t)|\mid
u-v\mid dxdt\geq 0. \eqno(5.8)
$$
Let $0<s<\tau<T$, and
$$
\eta (t)=\int_{\tau-t}^{s-t}\alpha_{\varepsilon}(\sigma)d\sigma ,\;\
\varepsilon <\min \{\tau, T-s\}.
$$
Here $\alpha_{\varepsilon}(t)$ is the kernel of mollifier with $
\alpha_{\varepsilon }(t)=0$ for $t\notin (-\varepsilon ,\varepsilon
)$. Then
$$
c\iint_{Q_T}|u-v|\eta(t)dxdt
+\int_{0}^{T}\left[\alpha_{\varepsilon}(t-s)-\alpha_{\varepsilon}(t-\tau)\right]|u-v|_{L^{1}(\Omega)}dt\geq
0.
$$
Let $\varepsilon \rightarrow 0$.  Then
$$
\int_{\Omega}|u(x,\tau)-v(x,\tau)|dx\leq
\int_{\Omega}|u(x,s)-v(x,s)|dx+c\int_{s}^{\tau}\int_{\Omega}|u-v|dxdt.\eqno(5.9)
$$
By the Gronwall inequality, we have
$$
\int_{\Omega}|u(x,\tau)-v(x,\tau)|dx\leq
c\int_{\Omega}|u(x,s)-v(x,s)|dx,
$$
 letting $s \rightarrow 0$, we have the conclusion.

 \textbf{Proof of Theorem 2.11} From proof of Theorem 2.10, we only need to deal with the term
 $$
\begin{gathered}
\int_{\Omega}\text{sgn}(u-v)(A(u,x,t)-A(v,x,t))\Delta\phi dx
\hfill \\
 \begin{array}{*{20}{c}}
{}
\end{array}
=-\eta(t)\int_{\Omega_{\lambda}}|A(u,x,t)-A(v,x,t)|\left[\frac{2}{\lambda^2}|\nabla d|^2+\frac{2(d-\lambda)}{\lambda^2}\Delta d\right] dx
\hfill \\
 \begin{array}{*{20}{c}}
{}
\end{array}
\leq -\eta(t)\int_{\Omega_{\lambda}}|A(u,x,t)-A(v,x,t)|\frac{2(d-\lambda)}{\lambda^2}\Delta d dx,
\hfill \\
\end{gathered}
$$
we have
 $$
\begin{gathered}
\lim_{\lambda\rightarrow 0}\left|\int_{\Omega_{\lambda}}|A(u,x,t)-A(v,x,t)|\frac{2(d-\lambda)}{\lambda^2}\Delta d dx\right|
\hfill \\
 \begin{array}{*{20}{c}}
{}
\end{array}
\leq c\lim_{\lambda\rightarrow 0}\frac{1}{\lambda}\int_{\Omega_{\lambda}}|A(u,x,t)-A(v,x,t)| dx
\hfill \\
 \begin{array}{*{20}{c}}
{}
\end{array}
=c\lim_{\lambda\rightarrow 0}\frac{1}{\lambda}\int_{\Omega_{\lambda}}\left|\int_{v}^{u}a(s,x,t)ds\right| dx
\hfill \\
 \begin{array}{*{20}{c}}
{}
\end{array}
=0.
\hfill\\
\end{gathered}\eqno(5.10)
$$
Then we have the conclusion.

\textbf{Proof of Theorem 2.12} Since we have imposed the partial boundary value condition
$$
u(x,t)=v(x,t)=0, (x,t)\in \Sigma_p\times(0,T),
$$
with
$$
\Sigma_p=\{x\in \partial \Omega: \sum_{i=1}^Nf_i(x)n_i<0\}\bigcup \{x\in \partial \Omega: \sum_{i=1}^Na_{x_i}(\cdot,x,t)n_i\neq 0\} \bigcup \{x\in \partial \Omega: a(\cdot,x,t)\neq 0\}, \eqno(2.5)
$$
From proof of Theorem 2.10-Theorem 2.11, we know (5.4)(5.5) and (5.10) are still true. We only need to deal with the term
$$
-f_i(x)\phi_{x_i}|u-v|
$$
in (4.12). In other words, since there is not the condition (2.10), the inequality (5.6) is not true. Actually, by the partial boundary value condition (2.12) with the expression (2.5), if we denote that
$$
\Omega_{1\lambda}=\{x\in \Omega_{\lambda}: f_i(x)d_{x_i}<0\},\eqno(5.11)
$$
then we have
$$
\begin{gathered}
-\lim_{\lambda\rightarrow 0}\int_{\Omega}f_{i}(x)\phi_{x_{i}}|u-v|dx
\hfill \\
 \begin{array}{*{20}{c}}
{}
\end{array}
=- 2\lim_{\lambda\rightarrow 0}\int_{\Omega_{\lambda}}f_{i}(x)\frac{(d-\lambda)d_{x_i}}{\lambda^2}\eta(t)|u-v|dx
\hfill \\
 \begin{array}{*{20}{c}}
{}
\end{array}
\leq- 2\lim_{\lambda\rightarrow 0}\int_{\Omega_{1\lambda}}f_{i}(x)\frac{(d-\lambda)d_{x_i}}{\lambda^2}\eta(t)|u-v|dx
\hfill \\
 \begin{array}{*{20}{c}}
{}
\end{array}
\leq -2\lim_{\lambda\rightarrow 0}\frac{1}{\lambda}\int_{\Omega_{1\lambda}}f_{i}(x)d_{x_i}\eta(t)|u-v|dx
\hfill \\
 \begin{array}{*{20}{c}}
{}
\end{array}
=2\eta(t)\int_{\Sigma_p}(-f_i(x)n_i)|u-v|d\sigma
\hfill \\
 \begin{array}{*{20}{c}}
{}
\end{array}
=0.
\hfill\\
\end{gathered}
\eqno(5.12)
$$
Similar as the proof of Theorem 2.11, we have the conclusion.
\bigskip

\section{The explanation of Definition 2.6}

Let us give a simple explanation of Definition 2.6 lastly.

Let $u_{\varepsilon}$ be the solution of the regularized equation
\begin{equation*}
\frac{\partial u}{\partial t}=\frac{\partial}{\partial x_i}\left(a(u,x,t)\frac{\partial u}{\partial x_i}\right)+\varepsilon \Delta
u+f_i(x)D_iu-c(x,t)u+g(x,t),\ \ (x,t)\in Q_{T},\eqno(6.1)
\end{equation*}
with the initial-boundary value conditions (3.2)-(3.3). Multiplying both sides of (6.1) by
$\varphi S_{\varepsilon}(u_{\varepsilon}-k)$
 and integrating it over $Q_{T}$ yields
$$
\begin{gathered}
\iint\nolimits_{Q_{T}}\frac{\partial u_{\varepsilon}}{\partial t}\varphi S_{\varepsilon}(u_{\varepsilon }-k)dxdt
\hfill \\
 \begin{array}{*{20}{c}}
{}
\end{array}
=\iint\nolimits_{Q_{T}}\frac{\partial}{\partial x_{i}}\left(a(u_{\varepsilon},x,t)\frac{\partial u_{\varepsilon}}{\partial x_{i}}\right)\varphi
S_{\varepsilon}(u_{\varepsilon}-k)dxdt
\hfill \\
 \begin{array}{*{20}{c}}
{}
\end{array}+\varepsilon \iint\nolimits_{Q_{T}}\Delta u_{\varepsilon}\varphi
S_{\varepsilon}(u_{\varepsilon}-k)dxdt
\hfill \\
 \begin{array}{*{20}{c}}
{}
\end{array}
+\iint\nolimits_{Q_{T}}\frac{\partial
(f_{i}(x)u_{\varepsilon})}{\partial x_{i}}\varphi S_{\varepsilon
}(u_{\varepsilon}-k)dxdt-\iint\nolimits_{Q_{T}}f_{ix_i}(x)u_{\varepsilon}\varphi S_{\varepsilon
}(u_{\varepsilon}-k)dxdt
\hfill \\
 \begin{array}{*{20}{c}}
{}
\end{array}
-\iint\nolimits_{Q_{T}}c(x,t)u_{\varepsilon}\varphi S_{\varepsilon}(u_{\varepsilon }-k)dxdt+\iint\nolimits_{Q_{T}}g(x,t)\varphi S_{\varepsilon}(u_{\varepsilon }-k)dxdt.
\hfill\\
\end{gathered}\eqno(6.2)
$$
Integration by parts, (6.2) gives
$$
\begin{gathered}
 \iint\nolimits_{Q_{T}}I_{\varepsilon}(u_{\varepsilon}-k)\varphi
_{t}dxdt+\iint\nolimits_{Q_{T}}A_{\varepsilon}(u_{\varepsilon},x,t,k)\triangle\varphi
dxdt  \hfill \\
 \begin{array}{*{20}{c}}
{}
\end{array}
-\iint\nolimits_{Q_{T}}f_i(x)I_{\varepsilon}(u_{\varepsilon}-k)\varphi
_{x_{i}}dxdt
-\varepsilon\iint\nolimits_{Q_{T}}\nabla
u_{\varepsilon}\cdot\nabla\varphi S_{\varepsilon}(u_{\varepsilon
}-k)dxdt \hfill \\
 \begin{array}{*{20}{c}}
{}
\end{array}
 -\varepsilon\iint\nolimits_{Q_{T}}\mid\nabla
u_{\varepsilon}\mid^{2}S_{\varepsilon}^{\prime}(u_{\varepsilon
}-k)\varphi dxdt
+\iint\nolimits_{Q_{T}}\int_{k}^{u_{\varepsilon}}a_{x_i}(s,x,t)S_{\varepsilon}(s-k)ds\varphi_{x_{i}}dxdt \hfill \\
 \begin{array}{*{20}{c}}
{}
\end{array}
-\iint\nolimits_{Q_{T}}a(u_{\varepsilon},x,t)\mid\nabla u_{\varepsilon}\mid^{2}h_{\varepsilon}(u_{\varepsilon}-k)\varphi dxdt \hfill \\
 \begin{array}{*{20}{c}}
{}
\end{array}
-\iint\nolimits_{Q_{T}}f_{ix_i}(x)(u_{\varepsilon}-k)\varphi S_{\varepsilon
}(u_{\varepsilon}-k)dxdt+\iint\nolimits_{Q_{T}}f_{ix_i}(x)\int_{k}^{u_{\varepsilon}}(s-k)h_{\varepsilon}(s-k)ds\varphi dxdt
\hfill \\
 \begin{array}{*{20}{c}}
{}
\end{array}
-\iint\nolimits_{Q_{T}}c(x,t)u_{\varepsilon}\varphi S_{\varepsilon}(u_{\varepsilon }-k)dxdt+\iint\nolimits_{Q_{T}}g(x,t)\varphi S_{\varepsilon}(u_{\varepsilon }-k)dxdt
\hfill \\
 \begin{array}{*{20}{c}}
{}
\end{array}
= 0.
\hfill\\
\end{gathered}\eqno(6.3)
$$

By discarding the terms
\begin{equation*}
-\iint\nolimits_{Q_{T}}a(u_{\varepsilon},x,t)\mid\nabla u_{\varepsilon}\mid^{2}S_{\varepsilon}^{\prime}(u_{\varepsilon}-k)\varphi dxdt,\eqno(6.4)
\end{equation*}
and
$$
-\varepsilon\iint\nolimits_{Q_{T}}\mid\nabla
u_{\varepsilon}\mid^{2}S_{\varepsilon}^{\prime}(u_{\varepsilon
}-k)\varphi dxdt
$$
in (6.3), we have
$$
\begin{gathered}
 \iint\nolimits_{Q_{T}}I_{\varepsilon}(u_{\varepsilon}-k)\varphi
_{t}dxdt+\iint\nolimits_{Q_{T}}A_{\varepsilon}(u_{\varepsilon},x,t,k)\triangle\varphi
dxdt  \hfill \\
 \begin{array}{*{20}{c}}
{}
\end{array}
-\iint\nolimits_{Q_{T}}f_i(x)I_{\varepsilon}(u_{\varepsilon}-k)\varphi
_{x_{i}}dxdt
-\varepsilon\iint\nolimits_{Q_{T}}\nabla
u_{\varepsilon}\cdot\nabla\varphi S_{\varepsilon}(u_{\varepsilon
}-k)dxdt \hfill \\
 \begin{array}{*{20}{c}}
{}
\end{array}
 +\iint\nolimits_{Q_{T}}\int_{k}^{u_{\varepsilon}}a_{x_i}(s,x,t)S_{\varepsilon}(s-k)ds\varphi_{x_{i}}dxdt \hfill \\
 \begin{array}{*{20}{c}}
{}
\end{array}
-\iint\nolimits_{Q_{T}}f_{ix_i}(x)(u_{\varepsilon}-k)\varphi S_{\varepsilon
}(u_{\varepsilon}-k)dxdt+\iint\nolimits_{Q_{T}}f_{ix_i}(x)\int_{k}^{u_{\varepsilon}}(s-k)h_{\varepsilon}(s-k)ds\varphi dxdt
\hfill \\
 \begin{array}{*{20}{c}}
{}
\end{array}
-\iint\nolimits_{Q_{T}}c(x,t)u_{\varepsilon}\varphi S_{\varepsilon}(u_{\varepsilon }-k)dxdt+\iint\nolimits_{Q_{T}}g(x,t)\varphi S_{\varepsilon}(u_{\varepsilon }-k)dxdt
\hfill \\
 \begin{array}{*{20}{c}}
{}
\end{array}
\geq0.
\hfill\\
\end{gathered}\eqno(6.5)
$$
Let $\varepsilon \rightarrow 0$. We can get
$$
\begin{gathered}
 \iint\nolimits_{Q_{T}}|u-k|\varphi
_{t}dxdt+\iint\nolimits_{Q_{T}}|A(u,x,t)-A(k,x,t)|\triangle\varphi
dxdt  \hfill \\
 \begin{array}{*{20}{c}}
{}
\end{array}
-\iint\nolimits_{Q_{T}}f_i(x)|u-k|\varphi
_{x_{i}}dxdt
\hfill \\
 \begin{array}{*{20}{c}}
{}
\end{array}
+\iint\nolimits_{Q_{T}}\int_{k}^{u}a_{x_i}(s,x,t)\text{sign}(s-k)ds\varphi_{x_{i}}dxdt \hfill \\
 \begin{array}{*{20}{c}}
{}
\end{array}
-\iint\nolimits_{Q_{T}}f_{ix_i}(x)(u-k) \text{sign}
(u-k)\varphi dxdt
\hfill \\
 \begin{array}{*{20}{c}}
{}
\end{array}
-\iint\nolimits_{Q_{T}}c(x,t)u\text{sign}(u-k)\varphi dxdt+\iint\nolimits_{Q_{T}}g(x,t) \text{sign}(u-k)\varphi dxdt
\hfill \\
 \begin{array}{*{20}{c}}
{}
\end{array}
\geq0.
\hfill\\
\end{gathered}\eqno(6.6)
$$

The inequality (6.6) is just the classical entropy inequality used in [3][5]etc. However,  the term (6.4) can not be thrown away casually. In fact, this term includes many information of the uniqueness [9-12], [18-19],[21-23][29]. The difficulty lies in that, when we let $\varepsilon \rightarrow 0$, what is the limit of the term (6.4) is very difficult to depict out, so it is almost impossible to remain the limit to the end, one has to throw it away [3][5].

In order to overcome this difficulty,  instead of multiplying both sides of (6.1) by
$\varphi S_{\varepsilon}(u_{\varepsilon}-k)$, we multiply both sides of (6.1) by
$\varphi S_{\eta}(u_{\varepsilon}-k)$, where $\eta$ is a small positive constant independent of $\varepsilon$. Then we can employ the weak convergent theory (Lemma 3.1), the uniqueness information of the term (6.4) remains, and we can prove the uniqueness of the entropy solutions by Kru\v{z}kov's method.

\section*{Availability of supporting data}
No applicable

\section*{Competing interests}
The author declares that he has no competing interests.

\bigskip
\section*{Funding} The paper is supported by Natural Science
Foundation of Fujian province (2019J01858),  and supported by SF of Xiamen University of Technology, China.

\section*{Author's contribution}

The author reads and approves the final manuscript.

\bigskip
\section*{Acknowledgement}
The author would like to think reviewers for their good comments.

 \bigskip

\bigskip
\section*{Reference}
\bigskip

\noindent[1] Oleinik  O.A. and Samokhin V. N., Mathematical Models in boundary Layer Theorem, Chapman and Hall/CRC, 1999.

\noindent[2] Wu Zhuoqun, Zhao Junning, Yin Jinxue and Li Huilai, Nonlinear Diffusion Equations, Word Scientific Publishing, Singapore, 2001.

\noindent[3] Vo\'{l}pert A. I. and Hudjaev, S. I., On the problem for quasilinear degenerate parabolic equations of second order (Russian), Mat.Sb.,  3(1967) 374-396.

\noindent[4] Zhao Junning, Uniqueness of solutions of quasilinear degenerate parabolic equations, Northeastern Math.J.,  1(2)(1985)153-165.

\noindent[5] Wu Zhuoqun and Zhao Junning, The first boundary value problem for quasilinear degenerate parabolic equations of second order in several variables, Chin.Ann. of Math.,  4B(1)(1983)57-76.

\noindent[6] Brezis H. and Crandall M.G., Uniqueness of solutions of the initial value problem for $u_{t}-\Delta \varphi (u)=0$, J. Math.Pures et Appl., 58(1979) 153-163.

\noindent[7] Kru\v{z}kov S. N., First order quasilinear equations in several independent variables, Math. USSR-Sb.,  10(1970)217-243.

\noindent[8] Cockburn B. and Gripenberg G., Continuous dependence on the nonlinearities of solutions of degenerate parabolic equations, J.  Differential Equations,  151(1999)231-251.

\noindent[9] Chen G. Q. and Perthame B., Well-Posedness for non-isotropic degenerate parabolic-hyperbolic equations, Ann. I. H. Poincare-AN,
20(2003)645-668.

\noindent[10] Chen G. Q. and DiBenedetto E., Stability of entropy solutions to Cauchy problem for a class of nonlinear hyperbolic-parabolic equations, SIAM J.Math.Anal., 33(4)(2001)751-762.

\noindent[11] Karlsen K.H. and Risebro N.H. On the uniqueness of entropy solutions of nonlinear degenerate parabolic equations with rough coefficient, Discrete Contain. Dye. Sys.,   9(5)(2003)1081-1104.

\noindent[12] Bendahamane M. and Karlsen, K.H., Reharmonized entropy solutions for quasilinear anisotropic degenerate parabolic equations, SIAM J. Math. Anal.,  36(2)(2004)405-422.

\noindent[13] Carrillo J., Entropy solutions for nonlinear degenerate problems, Arch. Rational Mech. Anal.,   147(1999) 269-361.

\noindent[14]  Li Yachun,  Wang Oin, Homogeneous Dirichlet problems for quasilinear anisotropic degenerate parabolic- hyperbolic equations, J.  Differential Equations,  252(2012)4719-4741.

\noindent[15] Lions P.L., Perthame B., and Tadmor E., A kinetic formation of multidimensional conservation laws and related equations, J. Amer. Math. Soc.,  7(1994) 169-191.

\noindent[16]  Kobayasi K.,  Ohwa H., Uniqueness and existence for anisotropic degenerate parabolic equations with boundary conditions on a bounded rectangle, J. Differential Equations,  252(2012)137-167.

\noindent[17] Yin J.,  Wang C., Evolutionary weighted $p-$Laplacian with boundary degeneracy, J. Differential Equations,  237 (2007) 421-445.

\noindent[18] Zhan Huashui, The Solutions of a Hyperbolic-parabolic mixed type equation on  half-space domain, J.  Differential Equations, 259(2015) 1449-1481.

\noindent[19] Zhan Huashui, The entropy solution of a hyperbolic-parabolic mixed type equation, SpringerPlus,  5(2016) 1811.

\noindent[20] Zhan Huashui, Reaction diffusion equations with boundary degeneracy, Electron. J. Differential Equations, 81(2016)1-13.

\noindent[21] Zhan Huashui, The study of the Cauchy problem of a second order quasilinear degenerate parabolic equation and the parallelism of a Riemannian manifold, Doctorial Dissertation, Xiamen University, 2004.

\noindent[22] Zhao Junning and Zhan Huashui, Uniqueness and stability of solution for Cauchy problem of degenerate quasilinear parabolic equations, Science in China Ser. A,  48(2005)583-593.

\noindent[23] Zhan Huashui,  On a hyperbolic-parabolic mixed type equation, Discrete and Continuous Dynamical systems, Series S., 10(3)(2017)605-624.

\noindent[24] Evans L.C., Weak convergence methods for nonlinear partial differential equations, Conference Board of the Mathematical Sciences, Regional Conferences Series in Mathematics Number 74,1998.

\noindent[25] Gu L. K., Second Order Parabolic Partial Differential Equations, Xiamen University Press, Xiamen, China, 2004.

\noindent[26] G. Fichera G., Sulle equazioni differenziali lineari ellittico-paraboliche del secondo ordine, Atti Accad, Naz. Lincei. Mem. CI. Sci Fis. Mat. Nat. Sez.  1(8)(1956)1-30.

\noindent[27] Oleinik  O. A.,  Radkevic E. V. , Second Order Differential Equations with Nonnegative Characteristic Form. Rhode Island: American Mathematical Society, and New York: Plenum Press, 1973.

\noindent[28]  Enrico G., Minimal Surfaces and Functions of Bounded Variation,
Birkhauser, Bosten. Basel. Stuttgart, Switzerland,1984.

\noindent[29] Zhan Huashui,  Feng Zhaosheng, Stability of hyperbolic-parabolic mixed type equations with partial boundary value condition, J.  Differential Equations, 264 (2018) 7384-7411.

\noindent[30]Abreu Eduardo, Colombeau Mathilde  and Panov Evgeny Yu, Approximation of entropy solutions to degenerate nonlinear parabolic equations, Z. Angew. Math. Phys. (2017) 68:133

\noindent[31] Frid Hermano, Li Yachun, A Boundary Value Problem for a Class of Anisotropic Degenerate Parabolic-Hyperbolic Equations, Arch. Rational Mech. Anal. 226 (2017) 975-1008.

\noindent[32]
Xin Z. and Zhang L., On the global existence of solutions to the
Pradtl's system, \emph{Adv. in Math.}, 191(2004), 88-133.

\noindent[33]
Amirat, Y., Bodart, O., Chechkin, G. A., Piatnitski, A. L., Boundary
homogenization in domains with randomly oscillating boundary.
\emph{Stoch. Proc. Appl.}, 121(2011),  1-23.

\noindent[34]
Basson, A., G¡äerard-Varet, D., Wall laws for fluid flows at a
boundary with random roughness. \emph{Comm. Pure Appl. Math.},
61(7)(2008),  941-987.

\noindent[35]
 Bucur, D., Feireisl, E., Nec¡¦asov¡äa, S.,
Wolf, J., On the asymptotic limit of the Navier-Stokes system on
domains with rough boundaries. \emph{J. Differ. Equ.},
244(11)(2008),  2890-2908.

\noindent[36]
Zhang J. W. and Zhao J. N., On the Global Existence and Uniqueness
of Solutions to Nonstationary Boundary Layer System, \emph{Science
in China, Ser. A}, 36(2006), 870-900.

\noindent[37]
Weinan E., Blow up of solutions of the unsteadly Pradtl's equations,
\emph{Comm. pure Appl. Math.}, 1(1997), 1287-1293.

\end{document}